\documentclass[11pt]{amsart}
\usepackage{latexsym,amssymb,amsmath}
\usepackage{epsfig}
\usepackage{epsf,float}
\newcount\cols
{\catcode`,=\active\catcode`|=\active
\gdef\Young(#1){\hbox{$\vcenter
{\mathcode`,="8000\mathcode`|="8000
\def,{\global\advance\cols by 1 &}%
\def|{\cr
      \multispan{\the\cols}\hrulefill\cr
       &\global\cols=2 }%
  \offinterlineskip\everycr{}\tabskip=0pt
  \dimen0=\ht\strutbox \advance\dimen0 by \dp\strutbox
    \halign
    {\vrule height \ht\strutbox depth \dp\strutbox##
      &&\hbox to \dimen0{\hss$##$\hss}\vrule\cr
     \noalign{\hrule}&\global\cols=2 #1\crcr
     \multispan{\the\cols}\hrulefill\cr%
   }
}$}}
}
\def\sqr#1#2{{\vcenter{\vbox{\hrule height.#2pt
\hbox{\vrule width.#2pt height#1pt \kern#1pt
\vrule width.#2pt}
\hrule height.#2pt}}}}

\def\s{\sigma}

\def\cf{\mathcal F}

\def\trank{\text{rank}}

\def\BC{\mathbb C}\def\BO{\mathbb O}\def\BS{\mathbb S}

\def\BP{\mathbb P}
\def\pp#1{\mathbb P^{#1}}
\def\fa{\mathfrak a}

\def\fd{\mathfrak d}

\def\pp#1{{\mathbb P}^{#1}}
\def\tdim{\rm dim}
\def\hd{,...,}
\def\ww{\wedge}
\def\upperp{{}^\perp}

\def\inv{{}^{-1}}

\def\cB{{\mathcal B}}
\def\cJ{{\mathcal J}}

\def\cF{{\mathcal F}}
\def\cQ{{\mathcal Q}}
\def\cR{{\mathcal R}}
\def\cS{{\mathcal S}}

\def\BB{{\mathbb B}}\def\BZ{{\mathbb Z}}

\def\11{\mathbf 1}

\def\FS{{\mathfrak S}}

\def\fe{{\mathfrak e}}

\def\ff{{\mathfrak f}}

\def\fg{{\mathfrak g}}
\def\fn{{\mathfrak n}}
\def\fp{{\mathfrak p}}

\def\ft{{\mathfrak t}}
\def\fl{{\mathfrak l}}
\def\l{\lambda}
\def\a{\alpha}

\def\o{\omega}

\def\s{\sigma}

\def\d{\delta}

\def\m{\mu}

\def\ot{{\mathord{\,\otimes }\,}}
\def\op{{\mathord{\,\oplus }\,}}

\def\lra{{\mathord{\;\longrightarrow\;}}}
\def\ra{{\mathord{\;\rightarrow\;}}}

\def\La#1{\Lambda^{#1}}

\def\tdim{\text{dim}\,}

\def\tmin{\text{ min }}

\def\trank{\text{rank}\,}

\newtheorem{theo}{Theorem}

\newtheorem{theorem}{Theorem}[section]
\newtheorem{proposition}[theorem]{Proposition}
\newtheorem{lemma}[theorem]{Lemma}

\theoremstyle{definition}

\newtheorem{example}[theo]{Example}

\theoremstyle{remark}

\begin{document}

\title{On secant   varieties  of Compact
Hermitian Symmetric Spaces}
\author{J.M. Landsberg
   and Jerzy Weyman}
\begin{abstract}{We show that the secant varieties
of   rank three compact Hermitian
symmetric spaces in their minimal homogeneous embeddings
are normal, with rational singularities. We show that their ideals
  are generated in degree three - with
one exception, the
secant variety of the $21$-dimensional spinor variety in $\pp{63}$ where
we show the ideal is generated in degree four.
We also discuss  the coordinate rings of secant
varieties of compact Hermitian
symmetric spaces.}
\end{abstract}
\thanks{Respectively supported by NSF grants DMS-0805782 and DMS-0600229}

\maketitle

\section{Introduction}
Let $K$ be an algebraically closed field of characteristic zero, let $V=K^{N+1}$
let $X\subset \BP V=\BP^N$ be a projective variety, and let
$\s(X)\subset \BP V$ denote its {\it secant variety}, the Zariski
closure of the set of points on the secant lines to $X$. Recently
there has been  interest in the ideals of
secant varieties of homogeneous varieties \cite{LMsec,LMsecb,LWsecseg,MR2319156, MR2310544,MR2252121,MR2201829,MR2168286}, and this paper contributes
to their study.

If the ideal of a variety $X$ is generated in degree two, the minimal possible
degree  of generators for the ideal of $\s(X)$ is three (\cite{LMsec} Cor. 3.2), although in general
one does not expect generators in degree three (e.g. this almost always fails
for complete intersections of quadrics). On the other hand, when $X$ is homogeneous, i.e.,
$V$ is an irreducible $G$-module where $G$ is a a semisimple algebraic group and
$X$ is the orbit of a highest weight line (so in particular, the ideal of
$X$ is generated in degree two), 
 in all previously known
examples (mostly just the rank two compact Hermitian symmetric spaces), the ideal
of $\s(X)$ was generated in degree three.  

In this paper we determine the generators of the ideals of the secant
varieties of rank three compact Hermitian symmetric spaces in their minimal 
homogeneous embeddings, which we abbreviate CHSS.
There is one suprise, the ideal of the secant variety of the $D_7$ spinor variety
is not generated in degree three, which answers a question posed in \cite{LMsec}, Section 3.
Recently, L. Manivel has made significant progress towards determining the generators
of the ideals of secant varieties of spinor varieties in general, see \cite{Manspin}.

While determining the generators of the ideals of secant varieties of higher rank 
CHSS seems out of reach at the moment, we show that for all other CHSS other
than spinor varieties, there are indeed generators in degree three. Moreover

\begin{theorem}\label{singthm} Let $X\subset \BP V$ be a rank three CHSS
in its minimal homogeneous embedding. Then $\s(X)$ is normal,
with rational singularities.   
\end{theorem}

Let $G/P_j\subset \BP V_{\o_j}$ denote the embedded rational homogeneous variety where $P_j$ is the maximal parabolic
associated to the simple root $\a_j$, using the ordering of the roots as in \cite{bour}.   Among the rank three CHSS are the Legendrian varieties,
$E_7/P_7, D_6/P_6, A_5/P_3=G(3,K^6), C_3/P_3=G_{Lag}(3,K^6), Seg(\pp 1\times Q)$ (where the last is the Segre product of a $\pp 1$ with a quadric hypersurface) which have the property that
their secant varieties are the ambient $\BP V$, so there is no need to
study their ideals and singularities. We let $S_{\pi}W$ denote the irreducible $SL(W)$ module associated to the partition $\pi$. The generators of the ideals 
in the remaining cases are as follows:

\begin{theorem}\label{g3nthm} Let $\tdim W\geq 7$.
The  ideal of the secant variety of the Grassmannian of $3$-planes in its Plucker embedding,  $\s(G(3,W^*))\subset \BP \La 3 W^*$
is   generated in degree three by the $SL(W)$-module of highest weight $2\o_1+\o_7$ occurring in $S^3(\La 3W)$, i.e., 
$S_{3,1^6}W=S_{3111111}W \subset S^3(\La 3 W) $.
\end{theorem}

Theorem \ref{g3nthm} is proved in \S\ref{g3sect}.

\begin{theorem}\label{spinthm} Let $X=\BS_7\subset \pp{63}$ be the $D_7$-spinor
variety. Then the ideal of $\s(X)$ is generated by
the irreducible $D_7$-module with highest weight ${\o_{4}}$ in degree four.  
\end{theorem}

Theorem \ref{spinthm} is discussed in \S\ref{spinsect}.

For a vector space $A$, let $K_S(A)\subset A^{\ot 3}$ denote the
kernel of the symmetrization map $S^2A\ot A\ra S^3A$, it is a $GL(A)$-module
isomorphic to $S_{21}A$. We let $\pi_S: (A\ot B)^{\ot 3} \ra (A\ot B)^{\ot 3}$
denote the symmetrization map whose image is $S^3(A\ot B)$.

\begin{theorem}\label{dcomthm} Let $X=Seg(\BP A^*\times Y)\subset \BP (A\ot W)^*$ be a reducible
rank three CHSS in its minimal homogeneous embedding (other than $Seg(\pp 1\times Q)$).
Let $I_3(\s(Y))\subset S^3W$ and $S_1(Y)\subset K_S(W)$ respectively denote the
modules generating the ideal of the secant variety of $Y$ and
the space of linear syzygies for the ideal of $Y$.
Then  the ideal of $\s(X)$ is generated in degree three
by 
$$\La 3 A\ot \La 3 W, \ S^3A\ot I_3(\s(Y)), \ 
{\rm and\ } \pi_S(K_S(A)\ot S_1(Y))
$$
Explicitly, the modules are
$$\begin{array}{c|c|c}
Y&  I_3(\s(Y)) &S_1(Y)\\
\hline   \\
G(2,6) & K & V_{\o_1+\o_5}^{A_n}\\
G(2,n+1), \ n\geq 6 &V_{\o_6}^{A_n} & V_{\o_1+\o_5}^{A_n}\\
Seg(\BP B^*\times \BP C^*) & \La 2B\ot \La 2 C & (S_{21}B\ot \La 2 C)\op (\La 3 B\ot S_{21}C)\\
  \BS_5 & 0 & V^{\fd_5}_{\o_4}\\
\BO\pp 2 & K & V^{\fe_6}_{\o_2}
\end{array}
$$
 
\end{theorem}

The degree three statement in Theorem \ref{dcomthm} is a consequence of the more general result:

\begin{proposition}\label{segcubicsc}
Let $Y\subset \BP A^*$ and $Z\subset \BP B^*$ be varieties and let 
$Seg(Y\times Z)\subset \BP (A^*\ot B^*)$ denote their Segre product.
Then 
\begin{align*}
&I_3(\s(Seg(Y\times Z))=\\
&
\La 3 A\ot\La 3 B\ \op \ \pi_S(S_1(Y)\ot K_S(B) \op K_S(A)\ot S_1(Z))\ \op\  I_3(\s(Y))\ot S^3B\ 
\op \ S^3A\ot I_3(\s(Z)).
\end{align*}
\end{proposition}

Theorem  \ref{dcomthm} and Proposition \ref{segcubicsc} are proven in \S\ref{dcomsect}.

For  higher rank irreducible CHSS we have the following result, which is proved
in \S\ref{gkvproof}:

\begin{proposition}\label{G(k,V)} The ideal of $\s(G(k,K^m))\subset \BP V_{\o_k}^*=\BP(\La k K^m)$ with $m\geq 2k$ contains
the irreducible $A_{m-1}$-modules $V_{2\o_{k-2}+\o_{k+4}}\op V_{ \o_{k-4}+2\o_{k+2}}\in S^3V_{\o_k}$ among its generators.
\end{proposition}

 In \cite{Manspin} it is shown that  $D_7/P_7, D_8/P_8$ are the only spinor varieties
 whose secant varieties have empty ideal in degree three, thus these are the only 
CHSS in not having cubics in the ideal of its secant variety. (Segre products
and Veronese re-embeddings of any varieties contain cubics in the ideals of
their secant varieties as there are cubics in the ideals of
the Segre products and Veronese embeddings of projective spaces.)

In \S\ref{coordringsect} we discuss the coordinate ring of $\s(X)$ for arbitrary
rational homogeneous varieties. A key point is that when $X\subset \BP V$ is homogeneous, 
$\s(X)$ is the closure of the orbit of
the sum of a highest weight vector and a lowest weight vector.

We obtain our results using the methods of \cite{weyman},
as described in Theorem \ref{weymanthm} below, along with some new
results about induced representations. In brief, in each
case we obtain a desingularization of $\s(X)$, by exploiting
that fact that each $X$ has a Legendrian \lq\lq smaller cousin\rq\rq ,
and apply Weyman's method to this desingularization.

 The case of $Seg(\BP A^*\times Q)\subset \BP (A^*\ot W^*)$ is immediate as
$\s (Seg(\BP A^*\times Q))=\s(Seg(\BP A^*\times \BP W^*))$.

\subsection*{Notation}\label{notation}
For a variety $Z\subset \BP V$, we let
$\hat Z\subset V$ denote the corresponding cone.
We adopt the following conventions: $K$ is an algebraically closed field
of characteristic zero, $G$ is a complex semi-simple algebraic group, $P$ a parabolic subgroup, $X=G/P\subset \BP V$ denotes a rational homogeneous
variety in its minimal homogeneous embedding.
We use German letters to denote Lie algebras
associated to algebraic groups.
We use the ordering of roots as in \cite{bour}.
   The fundamental weights
and  the simple roots  of $\fg$ are respectively
denoted  $\o_i$  and
$\a_i$.
  $P_k$ denotes the maximal parabolic of $G$ obtained by
deleting 
the root spaces corresponding to  negative   roots having a nonzero coefficient on the simple
root $\a_k$.
More generally,
for $J=(j_1\hd j_s)$,   $P_J$ denotes the parabolic obtained
by deleting the negative root spaces     having a nonzero coefficient on any of the simple
roots
$\alpha_{j_1}\hd \alpha_{j_s}$.
$\Lambda_{\fg},\Lambda_{G}$
respectively denote the weight lattices of $\fg$, $G$,
and $\Lambda^+_{\fg}\subset \Lambda_{\fg}$, $ \Lambda^+_G\subset \Lambda_G$ the dominant weights. We let $L\subset P$ be a (reductive)
Levi factor and $\ff=[\fl,\fl]$ a semi-simple Levi
factor.  We write $\fp=\fl +\fn$, where $\fn$
is nilpotent.  $V^{\fg}_{\l}$ denotes the irreducible $\fg$-module with
highest weight $\l$ and we often supress $\fg$ in the notation.
Unless otherwise noted, $G$ will be simply connected so $\Lambda^+_G=\Lambda^+_{\fg}$ and
we will freely switch from rational $G$-modules to $\fg$-modules.

When dealing with $\fa_n$-modules we sometimes use partitions
to index highest weights, with the dictionary
$\pi=(p_1\hd p_{n+1})$ corresponds to the
weight $(p_1-p_2)\o_1+(p_2-p_3)\o_2+\cdots + (p_n-p_{n+1})\o_n$.
We write $S_{\pi}K^n$ for the associated module.
Sometimes we abbreviate a partition
$(i_1\hd i_1,i_2\hd i_2\hd i_k\hd i_k)=((i_1)^{a_1},(i_2)^{a_2}\hd (i_k)^{a_k})$ where $i_s$ occurs $a_s$ times.

\subsection*{Acknowledgments} We thank W. Kra\'skiewicz and C. Robles for significant
help with computer calculations and L. Manivel for pointing out and correcting an
error in an earlier version of this paper.

\section{Method of proof}\label{sect2}

\subsection{The basic theorem of \cite{weyman}}

\begin{theorem}\cite{weyman}\label{weymanthm}
 Let $Y\subset \BP V$ be a
variety and suppose there is a
projective variety $\cB$ and a vector bundle
$q: E\ra \cB$ that is a subbundle of a trivial bundle
$\underline V \ra \cB$ with $\underline V_z\simeq V$ for
$z\in \cB$ such that $  \BP E\ra   Y$ is a desingularization
of $Y$.  Write $\eta=E^*$ and $\xi=(\underline V/E)^*$.

If
the sheaf cohomology groups
$H^i(\cB,S^d\eta)$ are all zero for $i>0$ and $d>0$
and if the linear maps
$H^0(\cB,S^d\eta)\ot V^*\ra H^0(\cB,S^{d+1}\eta)$
are surjective for all $d\geq 0$, then
\begin{enumerate}
\item
$\hat Y$ is normal, with rational singularities   

\item The coordinate ring $K[\hat Y]$ satisfies
  $K[\hat Y]_d\simeq H^0(\cB,S^d\eta)$.

\item The 
vector space of minimal generators of the ideal of $Y$ in
degree $d$ is isomorphic to  
$H^{d-1}(\cB,\La{d }\xi)$ which is also the homology of the
complex
\begin{equation}\label{complex}
0\ra\La 2V \ot H^0(\cB, S^{d-2}\eta)   \lra   V \ot H^0(\cB, S^{d-1}\eta)    \lra  
 H^0(\cB, S^d\eta)\ra 0.
\end{equation}

\item More generally,  
$\oplus_j H^j(\cB, \La{i+j}\xi)$ is isomorphic to
the $i$-th term in the minimal free resolution
of $Y$.

\item If moreover $Y$ is a $G$-variety
and the desingularization is $G$-equivariant,
then the identifications above are as $G$-modules.
\end{enumerate}
\end{theorem}

\subsection{The basic theorem applies in our case}
In this paper the desingularizations will all be
by a homogeneous bundle $\BP E$ such that the corresponding
bundle $\eta$ is irreducible. In this case 
we have:
\begin{proposition}\label{acyclicprop} Notations as above, if
$Y$ is a $G$ variety, $\cB=G/P$ and $\eta$ is
induced from an irreducible $P$-module, then
the sheaf cohomology groups
$H^i(\cB,S^d\eta)$ are all zero for $i>0$
and   the linear maps
$H^0(\cB,S^d\eta)\ot V^*\ra H^0(\cB,S^{d+1}\eta)$
are surjective for all $d\geq 0$. In particular
  all the conclusions of \ref{weymanthm} apply.
\end{proposition}

\begin{proof} An irreducible homogeneous bundle can have nonzero
cohomology in at most one degree, but a quotient
bundle of a trivial bundle has nonzero sections,
thus $H^0(B, \eta)$ is a nonzero irreducible module
and all other $H^j(B,\eta)$ are zero. 
Let $\ff\subset\fp\subset \fg$ be a semi-simple
Levi factor, so the weight lattice of $\ff$ is a sublattice
of the weight lattice of $\fg$,   let $\ft^c$ denote
the complement of $\ft_{\ff}$ (the torus of $\ff$) in $\ft_{\fg}$ and  
let $\fg_0=\ff + \ft^c$ denote the Levi factor of $\fp$.
   $\eta$ is induced from an
irreducible  $\fg_0$-module $U$ which is a weight space
for $\ft^c$ having 
non-negative   weight, say $(w_1\hd w_p)$.
  The bundle
  $\eta^{\ot d}$,  corresponds to a module
which is $U^{\ot d}$ as an $\ff$-module and 
  is a weight space with weight $(dw_1\hd dw_p)$ for
the action of $\ft^c$.
Thus $S^d\eta$ is completely reducible and each component of $S^d\eta$ is
very ample and in particular acyclic.

To prove the second assertion, 
  consider the maps 
$V^* \ot H^0(B,S^{r-1}\eta)\ra H^0(B,S^{r }\eta)$.
Note that $H^0(B,S^{j}\eta)\subset S^{j}V^*$.
The proof of Proposition \ref{acyclicprop}
will be completed by Lemma \ref{littlem} below
applied to $U$   and
each irreducible component of $H^0(B,S^{r}\eta)$.\end{proof}

Let $M^{\fg}_{\fg_0}$ denote the sub-category  of the category of $\fg_0$-modules  generated 
under direct sum by
the irreducible $\fg_0$ modules with highest weight in $\Lambda^{+}_{\fg}\subset \Lambda^+_{\fg_0}$
and note that it is closed under tensor product. 
Let $M_{\fg}$ denote the category of $\fg$-modules.
Define an additive functor $\cf : M^{\fg}_{\fg_0}\ra M_{\fg}$ which takes an irreducible
$\fg_0$-module with highest weight $\l$ to the corresponding irreducible
$\fg$-module with highest weight $\l$.

\begin{lemma}\label{littlem} Let $\fl\subset \fg$ and $\cF$ be as
above. Let $U ,W$ be irreducible $\fl$-modules Then
$$\cF(U \otimes W)\subseteq \cF(U )\ot \cF(W).$$
\end{lemma}
\begin{proof}  Let $N\subset P$ denote the unipotent radical of $P$. Any $L$-module $W$ may be considered
as a $P$-module where $N$ acts trivially. Saying $V=\cf(W)$ means that $V$ is the $G$-module
parabolically induced from $W$ and $W$ is the set of $N$-invariants of $V$.  The $N$-invariants of $\cf(U)\ot \cf(W)$ contain $U\ot W$. \end{proof}

\section{Proof of Proposition \ref{G(k,V)}}\label{gkvproof}

 We recall some facts from \cite{LMsec}, Section 3.
For any variety $X\subset \BP W^*$ whose ideal
is generated in degree two,  $I_3(\s(X))=S^3W\cap (I_2(X)\ot W)$
with the intersection being taken inside $S^2W\ot W$.

Let $G$ be semi-simple, let  $W=V_\l$ and let $X=G/P\subset\BP W^*$ be the orbit of a highest weight line.
Assume that an irreducible $G$-module $V_\m$
appears in $S^3V_\l$   and does not
appear in $V_{2\l}\ot V_\l$, where $V_{2\l}$ is identified with the Cartan square of $V_\l$ (i.e.,
the unique submodule of $S^2 V_\l$ isomorphic to $V_{2\l}$). Then $V_\m$  must be in $I_3(\s(X))$
as $S^2V_\l =I_2(X)\op V_{2\l}$. 

\medskip

Let $G=SL(n,K)$ and $W=\La kK^m=V_{\o_k}$.
It follows from the Pieri formulas   
that the modules $V_{2\o_{k-2}+\o_{k+4}}$ and $V_{ \o_{k-4}+2\o_{k+2}}$  do not occur in $V_{2\o_k}\ot V_{\o_k}$. 
To see that they occur in $S^3W$, 
identify $V_{\o_k}$ with $\La k (K^n)$ where $K^n$ has a basis $\lbrace e_1 ,\ldots ,e_n\rbrace$.
First observe that $\La 6 (K^6 )\subset S^3 (\La 2(K^6 ) )$, in fact if
$(f_1\hd f_6)$ is a basis of $K^6$, then the inclusion takes it to the Pfaffian
$$f_1\ww\cdots \ww f_6\mapsto
\sum_{\sigma} {\rm sgn}(\sigma )(f_{\sigma (1)}\ww f_{\sigma (2)})\circ (f_{\sigma (3)}\ww f_{\sigma (4)})\circ (f_{\sigma (5)}\ww f_{\sigma (6)})
$$
where we sum over all permutations $\sigma\in\FS_6$ satisfying
$$\sigma (1)<\sigma (2) ,\  \sigma (3)<\sigma (4) ,\ \sigma (5)<\sigma (6),\  \sigma (1)<\sigma (3) <\sigma (5).$$

Now considering $K^6\subset K^n$
as the span of $\lbrace e_{k-1},\ldots ,e_{k+4}\rbrace$ we can produce a highest weight vector of
$V_{2\o_{k-2}+\o_{k+4}}$ by wedging each term $f_i\ww f_j=e_{i+k-2}\ww e_{i+k-2}$ in the summation
with $e_1\wedge\ldots\wedge e_{k-2}$. We leave it to the reader to check
the resulting vector has the desired properties.
The module $V_{\o_{k-4}+2\o_{k+2}}$ occurs in $S^3 W$ as well by symmetry (or one can
define an analogous map).

\section{Desingularizations for secant varieties 
of Rank 3 CHSS}\label{desingsect}
In our situation the desingularizations are
  based on the observation
that in each case $X$ is swept out by the union
of Legendrian varieties $X_{small}$ and 
$\s(X)$ is the union of the $\s(X_{small})$'s which
are linear spaces.

  Here is a table of $X,X_{small},\cB$ and the desingularizing bundle 
$E$ over
$B$:
$$\begin{array}{c|c|c|c}
X\ \ & X_{small}&\cB={\rm base}\ \ \ &E={\rm bundle}\\
\hline \\
G(3,U) &G(3,6) &G(6,U)
& \La 3\cR_U
\\
\BS_7 &\BS_6 &Q^{12}
&  Spin(p\upperp /\hat p)
\\
\BP A\times Y &\pp 1\times Q&G(2,A)\times Y_Q 
&\cR_A\ot S  
 \end{array}
 $$

Here, if $Y=G/P\subset \BP W$, $Y_Q=G/P'$ is the variety obtained via
Tits' {\it shadows} that parametrizes a space of quadric
sections of $Y$ (see \cite{LMmagic}, \S 2). One takes the marked
Dynkin diagram for $Y\subset \BP W$ and looks for the largest subdiagram
whose resulting marked diagram is  a quadric hypersurface.
The marked diagram of $Y_Q$ is obtained by marking all nodes
bounding the nodes of the subdiagram corresponding to the quadric
hypersurface.
The irreducible homogeneous bundle $S\ra Y_Q$ is obtained from
the
irreducible  $\fp$-module of highest weight equal to the highest
weight of $W$ (where we have included the weight lattice
of the Levi factor of $\fp$ into the weight lattice of $\fg$).
  As such, it
is an irreducible sub-bundle of the trivial bundle with fiber $W$
and satisfies the hypotheses of Proposition \ref{acyclicprop}.

\begin{example}
 \begin{figure}[H]
\begin{center}
\epsfxsize=3in\epsfbox{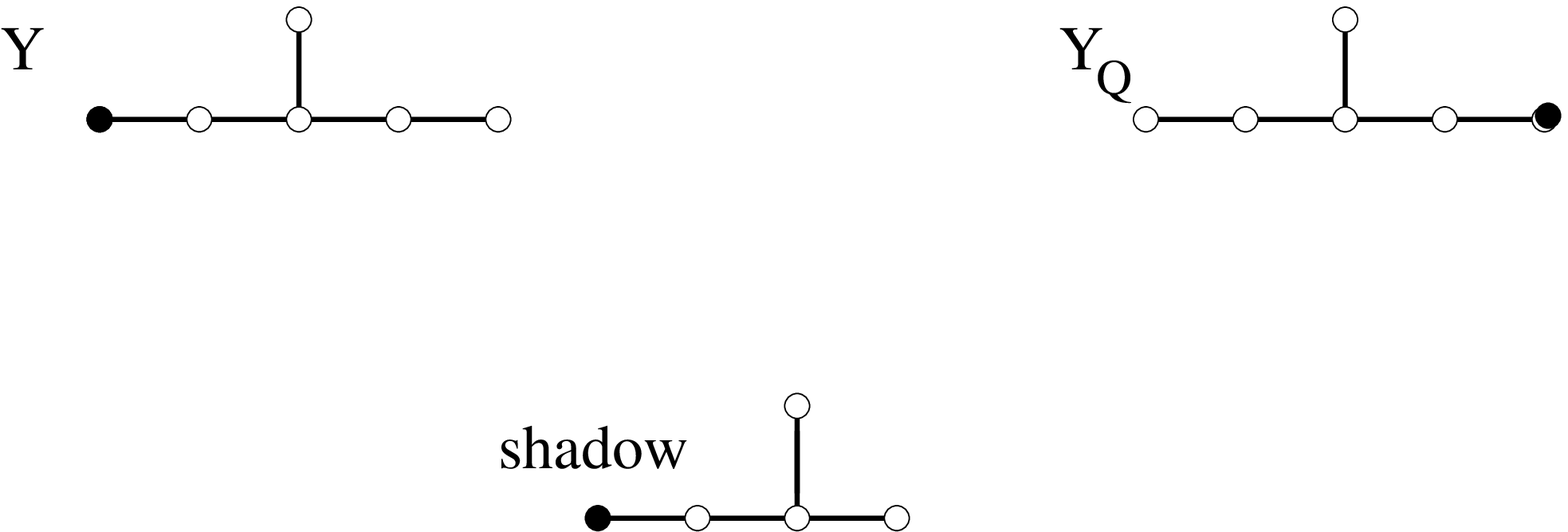}
  \caption{Shadow of $Y_Q=E_6/P_6$ on $Y=E_6/P_1$ is a $Q^8=D_5/P_1$}
\end{center}
\end{figure}
\end{example}

We have
$$
\begin{array}{c|c|c}
Y & Y_Q &\tdim Q\\
\hline \\
\BP B\times \BP C & G(2,B)\times G(2,C) & 2\\
G(2,W) & G(4,W) & 4\\
\BS_5 & Q^8 & 6\\
E_6/P_1 & E_6/P_6& 8
\end{array}
$$

When $\tdim Q=2^k$, we have a uniform model, $W=  \cJ_n(\BB)$,
where $\cJ_n(\BB)$ is the Jordan algebra of $n\times n$ $\BB$-Hermitian symmetric matrices.
In this model $\s (Y)$ is the set of
rank at most $2$ elements  with ideal generated
  by the $3\times 3$ minors.
(In the case $\BB=\BO$ the octonions, 
we have $n=3$ and care must be
taken when defining the determinant.)
In these cases the fiber of $S$ is isomorphic to
$\cJ_2(\BB)$.

When $Y=\BS_5$ we have $W\simeq \cS_5\simeq \La{even}K^5$
and the fiber of $S$ isomorphic to $\La{even}K^4$.

\medskip

\begin{lemma}\label{lemma31} Let $\tilde \s$ denote the total space of $E$. The image of $q:\tilde\s\ra V$ is $\hat \s(X)$
and $q:\tilde \s\ra \hat \s(X)$ is a resolution of singularities of $\hat \s(X)$.
\end{lemma}

\begin{proof}
In each case, the fiber $E_x\subset V_x=V$ over $x\in \cB$ is $\hat \s( X_{small,x})\subset E_x\subset V$
and $X_{small,x}\subset X$ by construction, so the image
of $q$ is contained in $\hat \s (X)$. On the other hand,
they both have the same dimension and $\hat \s (X)$ is 
reduced and irreducible. Thus we need only show the map
is generically one to one. It is clear that restricted
to each fiber $q$ is generically one to one, so it
is sufficient to show that there is a unique fiber over
a general point. For $\s(
G(3,W))$, a general point determines a unique
$6$ plane in $W$. For $\s(Seg(\BP A\times Y))$,
 a general point
  clearly determines a unique $2$-plane $A'$ in $A$. 
When $Y\neq \BS_5$, two elements
of $Y$ in general position lie in a unique $\cJ_2(\BB)$ 
and the unique fiber is $A'\ot \cJ_2(\BB)$. 
(The intersection   $Y\cap \BP \cJ_2(\BB)$ is the set of rank one elements in $\cJ_2(\BB)$,
i.e., a quadric of dimension $\tdim \BB$, see,
e.g., \cite{zak}, chapter VI or \cite{chaput}.)
For the case $Y=\BS_5$, fix  an
isotropic line
$L\subset K^{10}$, the set $\{ F\in \BS_{5}\mid L\subset F\}$
is the shadow of $L$ in $\BS_{5}$   and the span of this shadow
is the image of the fiber over $[L]\in Q^8$ in $K^{16}$.
Two general points $F,F'\in \BS_5$,
considered as $\pp 4$'s in $Q^8$ will intersect in
a point, i.e., an isotropic line in $K^{10}$, which
determines the unique fiber above a 
general point in their linear span.  

Finally for the case $X=\BS_7$, the same argument
for $\BS_5$ applies, as it is still true two general
$F,F'\in \BS_7$ will intersect in an isotropic line.
\end{proof}

Lemma \ref{lemma31} combined with  Theorem \ref{weymanthm}
and Proposition \ref{acyclicprop}
prove Theorem \ref{singthm}.

We now proceed with a case by case study.

\section{Case of $X=G(3,W)$}\label{g3sect}

Let 
$$
R_p(\La k W^*)=\{ T\in \La k W^* \mid \exists K^p\subset W^*
{\rm{\ such \ that\ }} T\in \La kK^p\}
$$
$R_p(\La k W^*)$ is called a
  {\it rank variety} (or {\it subspace variety}). Such varieties are discussed in detail in (\cite{weyman},
\S 7). Their ideals are easy to describe, namely
$I_d(R_p(\La k W^*))$ consists of all modules  corresponding to copies of $S_{\pi}W$
occurring in $S^d(\La k W)$ where $\ell (\pi)>p$. However it is in general difficult to determine
generators of the ideal. $R_p(\La k W^*)$ is desingularized by
$\La k\cS\ra G(p,W^*)$, and the corresponding bundle $\xi=(\La k\underline{V}/\La k \cS)^*$
in general is not irreducible. When $p=\tdim W-1$ however $\xi$ is irreducible, which will be the key
to our proof.
 
 \begin{proposition}   $\hat\s(G(3,W^*))=R_6(\La 3 W^*)$. 
 \end{proposition}
 \begin{proof}
 A general point of $\hat\s(G(3,W^*))$ is of the form $v_1\ww v_2\ww v_3+v_4\ww v_5\ww v_6$, so
   $\hat\s(G(3,W^*))\subseteq R_6(\La 3 W^*)$ because the latter is compact and the former connected.
 But both varieties are of the same dimension $6(\tdim W)-17$ and are reduced and irreducible so they must be equal.
 \end{proof}

Write $V=\La 3W$.
We have 
$$S^3V= S_{3^3}W\op S_{32^21^2}W\op S_{2^21^5}W\op S_{31^6}W
$$ 
The last module corresponds to   a partition of length seven  and thus by the remark above, it is among the generators
of $I(R_6(\La 3 W^*))$ (because the ideal in degree two of any secant variety is empty), and the rest are not as their partitions
  have length at most six.

To show there are no generators in degree greater than three,
we need to prove
exactness in the middle step of \eqref{complex}  which in this case is:
$$ 
(S_{1^6}W \op S_{2211} W )\ot   
S^{r-2}(S_{111} W  )\mid_{|\pi |\leq 6}   \ra    
S_{111} W  \ot  S^{r-1}(S_{111} W )\mid_{|\pi |\leq 6}   \ra   
 S^{r }(S_{111} W )\mid_{|\pi |\leq 6} 
$$
for $r>3$.
 
The largest partition that can show up in
the middle has length nine, so once we have solved the problem
for $G(3,K^9)$ we are done.

Thus one could proceed calculate $H^d(G(6,K^9),\La{d+1}\xi)$ with
the aid of a computer to conclude (although the passage from
the cohomology of $\La{d+1}gr(\xi)$ to $\La{d+1}\xi$ might require
some effort). We will proceed differently, resolving the cases
of $\tdim W=7,8,9$ iteratively using   rank varieties with $p=\tdim W-1$.

For $\tdim W=7$, the method in \cite{weyman}, \S 7.3 shows that the ideal of
$R_6(\La 3 W)$ is generated by $S_{31^6}W$ and we are done.
For the next two cases  we proceed indirectly, calculating the
ideal of $R_7(\La 3K^8)$ (resp. $R_8(\La 3K^9)$), 
and show
these are in the ideal generated by $S_{31^6}W$
to complete the proof.

\begin{proposition}\label{5.1}
The  ideal of the rank variety $R_6(\La 3K^{7*})$
is generated in degree three by $S_{3,1^6}K^{7}$
included in $S^3(\La 3 K^{7})$ as described in the recipe in the proof.

The  ideal of the rank variety $R_7(\La 3K^{8*})$
is generated in degree five by $S_{3,2^2,1^5}K^8$
included in $ S^5(\La 3 K^8)$ as described in the recipe in the proof.

The  ideal of the rank variety $R_8(\La 3K^{9*})$
is generated in degrees four and
five  by $S_{4,1^8}K^9$   and $S_{3^2,2^2,1^5}K^9$
respectively
included in $S^4(\La 3 K^9)$
and $S^5(\La 3 K^9)$ as described in the recipe in the proof.
 \end{proposition}

\begin{proof} Thanks to the irreducibility of $\xi$ and its exterior powers, determination of modules generating the ideal
is a straightforward application of  the methods of \cite{weyman} and
is left to the reader. It remains
to show the above modules are all in the ideal generated by $S_{31^6}W$. To do this
we give explicit descriptions of the modules as spaces of polynomials.

We will encode the representations occurring in the $d$-th symmetric powers of $\La 3 W$ by Young tableaux of shape $\pi =\lambda(D)$ with $3d$ boxes, filled with the numbers $1,\ldots ,d$ with each number occurring three times. These tableaux are also assumed  to be weakly increasing in rows and strictly increasing in columns. We associate to such tableau $D$ the map $\rho(D)$ 

$$\rho (D) :\La{d_1'}W\otimes\cdots \otimes \La{d_r'}W \ra S_d(\La 3 W)$$
where $\pi'=(d_1'\hd d_r')$ is the conjugate partition to $\pi$.

The map $\rho(D)$  is defined as the composition of the following maps:
\item{a)} assuming there are $e_{i,s}$ boxes filled with $s$ in the $i$-th row, apply the embedding
$$\La{d_i'}W\ra \La{e_{i,1}}W\ot\ldots\ot\La{e_{i,d}}W$$
for each row of $D$,
\item{b)} Noting that for each $s$, 
$e_{1,s}+\cdots + e_{r,s}=3$, wedge the factors coming from different rows corresponding to the same number $s$ in $D$, i.e., after rearranging the factors define the projection to $(\La 3 W)^{\ot d}$ by sending, for each $s$,
$$\La{e_{1,s}}W\ot\ldots\ot\La{e_{r,s}}W\ra\La{3}W$$
 and tensoring the results.
\item{c)} Project $(\La 3W)^{\ot d}\ra S^d (\La 3 W)$ by symmetrizing.

We call the Young diagram $D$ {\it  the numbering scheme} associated to the map $\rho(D)$.
Write $\lambda(D)$ for the weight of the numbering scheme whose $i$-th entry is the length of the $i$-th column of the Young diagram of $D$ - this will be the highest weight of the associated module.

The four Schur functors mentioned in the statement of the lemma correspond to four numbering schemes.

$$D_1=\Young (1,1,1,2,2,3,3|2|3).$$

$$D_2=\Young (1,1,1,2,2,3,3,4|2,4,4|3).$$

$$D_3=\Young (1,1,1,2,2,3,3,4,4|2|3|4).$$

$$D_4=\Young (1,1,1,2,2,3,3,4,5|2,4,4,5|3,5).$$

It is clear that the images of the corresponding maps are in the ideals of corresponding rank varieties because of the length of the first row of each numbering scheme.

Decomposing the domain and range of $\rho(D_i)$ into irreducible representations, in all four cases    $S_{\lambda(D_i)}W$   is the only Schur functor occurring in both the domain and range of $\rho(D_i)$, and it occurs there with multiplicity one.
Thus
it only remains to see the maps $\rho(D_j)$ are nonzero, which is
the purpose of the following lemma.

 \end{proof}

\begin{lemma}\label{42} The numbering schemes $D_1, D_2, D_3, D_4$ all yield nonzero
modules. 
\end{lemma}

\begin{proof}
We prove that the image of the map
$$\rho(D_2 ): \La 8 W\otimes\La 3 W\otimes W\rightarrow S^4 (\La 3 W)$$
corresponding to the numbering scheme
$$D_2=\Young (1,1,1,2,2,3,3,4|2,4,4|3)$$
is nonzero in $S_4 (\La 3 W)$. The other cases are similar.

 Consider the contribution to the monomial $(e_1\wedge e_2\wedge e_3)(e_1\wedge e_2\wedge e_3)(e_1\wedge e_4\wedge e_5)(e_6\wedge e_7\wedge e_8)$ in the image of highest weight vector $\rho(D_2 )(e_1\wedge e_2\wedge e_3\wedge e_4\wedge e_5\wedge e_6\wedge e_7\wedge e_8)\ot (e_1\wedge e_2\wedge e_3)\ot (e_1)$. All occurrences of this monomial can be divided to 24 classes (corresponding to permutations of $\lbrace 1,2,3,4\rbrace$) according to the order in which the factors appear in $(\La 3 W)^{\ot d}$ after applying parts a) and b) of
the  definition of $\rho(D_2 )$.
In fact only two classes out of 24 are non-empty.
The factor $e_6\wedge e_7\wedge e_8$ has to come from the first factor, and one of the factors $e_1\wedge e_2\wedge e_3$ has to come from the fourth factor. The factor $e_1\wedge e_4\wedge e_5$ can come from the third factor (and this gives contribution $3$ to the coefficient) or from the second factor (and this gives contribution $1$ to the coefficient). Thus the coefficient is nonzero and and therefore $\rho(D_2 ) \neq 0$.

\end{proof}

To finish the proof of Theorem \ref{g3nthm}  we need to show that the ideal generated by the first
module contains the other modules. But this is clear by the definition of maps $\rho(D)$ and by 
the last part of the proof of Proposition \ref{5.1}   as the other numbering schemes  all contain
 the first.

\section{Case of $X=\BS_7$}\label{spinsect}
A desingularization of $\hat\s(\BS_7)
\subset V_{\o_7}^{D_7}$ is given
by the sub-bundle of  $q:E\ra Q^{12}=D_7/P_1$ whose fiber is isomorphic 
to $\hat\s(\BS_{6})=V_{\o_{6}}^{D_{6}}$. 
Thus we can apply   Theorem \ref{weymanthm}.

Note that $V_{\o_7}^{D_7}$ decomposes
to  $V_{\o_{6}}^{D_{6}}\op V_{\o_{5}}^{D_{6}}$ as a $D_6$ module, and this splitting gives
rise to the 
   the bundles $\xi$ and $\eta$ over $Q^{12}$. Thus they are
both irreducible and dual to one another.
In this case it is straightforward to calculate
$H^j(\La{j+1}\xi)$  if one knows the decomposition
of $\La{j+1}V_{\o_{6}}^{D_{6}}$. In fact
 we calculated the entire minimal free resolution
which is available at http://www.math.neu.edu/$\sim$weyman/mathindex.html for the interested
reader. In particular the only generator of the ideal
is the module $V_{\o_4}$ as stated in the theorem.

\section{Case   of $\s(Seg(\BP A^*\times Y))$}\label{dcomsect}

\begin{proof}[proof of Proposition \ref{segcubicsc}]
All spaces discussed in this section are to be considered as linear subspaces of $(A\ot B)^{\ot 3}=A^{\ot 3}\ot B^{\ot 3}$
and all evaluations are as multi-linear forms.
In particular, the symmetrization map 
\begin{align*}
\pi_S: (A\ot B)^{\ot 3}&\ra (A\ot B)^{\ot 3}\\
a_1\ot a_2\ot a_3\ot b_1\ot b_2\ot b_3&\mapsto \frac 16\sum_{\s\in \FS_3}a_{\s(1)}\ot a_{\s(2)}\ot a_{\s(3)}
\ot b_{\s(1)}\ot b_{\s(2)}\ot b_{\s(3)}
\end{align*}   realizes $S^3(A\ot B)\subset (A\ot B)^{\ot 3}$
as $\pi_S((A\ot B)^{\ot 3})$. Similarly we regard $S^2A\ot A\subset A^{\ot 3}$ as the image
of the symmetrization map $x\ot y\ot z\mapsto \frac 12(x\ot y\ot z+y\ot z\ot z)$ and likewise
for $S^2B\ot B\subset B^{\ot 3}$.

As mentioned in \S\ref{gkvproof}, for any variety $X\subset \BP V$, $I_3(\s(X))=(I_2(X)\ot V^*)\cap S^3V^*$.
Here we have the $GL(A)\times GL(B)$ decomposition of 
$S^3(A\ot B)\subset (A\ot B)^{\ot 3}$:  
$$S^3(A\ot B)=\La 3 A\ot \La 3 B\op\pi_S( K_S(A)\ot K_S(B))
\op S^3A\ot S^3B,
$$
 where $K_S(A)$ is the kernel  of
the map  $S^2A\ot A\ra S^3A$, which is a $GL(A)$-module isomorphic to $S_{21}A$.
In particular, if $R\in K_S(A)$, we have $R(u,v,w)=R(v,u,w)$ for all $u,v,w\in A^*$ and
\begin{equation}\label{keyproperty}
R(u,v,u)= -\frac 12 R(u,u,v) \ \ \forall u,v,\in A^*
\end{equation}
which holds because $R(u,v,w)+R(u,w,v)+R(v,w,u)+R(v,w,u)+R(w,u,v)+R(w,v,u)=0$, and setting
$w=u$ gives $2R(u,v,u)+2R(u,u,v)+2R(v,u,u)=4R(u,v,u)+2R(u,u,v)=0$.

\smallskip

The   factor $\La 3 A\ot \La 3 B$ is in $I_3(\s(Seg(Y\times Z))$ because it is 
$I_3(\s(Seg(\BP A^*\times\BP B^*))$.

\smallskip

By polarization and symmetry, a polynomial $P\in S^3V^*\subset V^{*\ot 3}=\{ {\rm{trilinear\ maps\ }} V\times V\times V\ra \BC\}$ is in the ideal of $\s(X)$ if and only if
$P(x_1,x_1,x_2)=0$ for all $x_1,x_2\in\hat X$. In what follows we consider
$P$ as a trilinear form on $V$.

\smallskip

By definition   $S_1(Y)=\{ T\in I_2(Y)\ot A \mid \pi_{S,A}(T)=0\}$, where 
$$\pi_{S,A}(a_1\ot a_2\ot a_3)= \sum_{\s\in \FS_3}a_{\s(1)}\ot a_{\s(2)}\ot a_{\s(3)}.
$$
 Elements of   
$S_1(Y)\subset K_S(A)$  have the property that as   tri-linear
forms they vanish on any triple of the form $(v,v,w)$ with $[v]\in Y$ and $w$ arbitrary.

 Any element of $K_S(A)\ot K_S(B)\subset (A\ot B)^{\ot 3}$ 
may be written as a sum $\sum R_i\ot T_i$ with $  R_i\in K_S(A)$ and $  T_i\in    K_S(B)$, thus
any element of $\pi_S(K_S(A)\ot K_S(B))$ is of the form
  $P=\pi_S(\sum R_i\ot T_i)$ with $  R_i\in K_S(A)$ and $  T_i\in    K_S(B)$.
We compute
\begin{align*}
\frac 12
P(v\ot z,v\ot z, w\ot y)&= \sum R_i(v,v,w)T_i(z,z,y)+ \sum R_i(v,w,v)T_i(z,y,z)+ \sum R_i(w,v,v)T_i(y,z,z)\\
&
=\sum R_i(v,v,w)T_i(z,z,y)+ 2\sum R_i(v,w,v)T_i(z,y,z)\\
&
=\sum R_i(v,v,w)T_i(z,z,y)+ 2\sum (-\frac 12 R_i(v,v,w)(-\frac 12 T_i(z,z,y))\\
&
=\frac 32\sum R_i(v,v,w)T_i(z,z,y).
\end{align*}
The first equality holds because of the six permutations in $\FS_3$, only three yield different elements,
  the second because $R_i\in S^2A\ot A$ and $T_i\in S^2B\ot B$, and the third by \eqref{keyproperty}.

First we show that $S_1(Y)\otimes K_S(B)\subset I_3(\s(Seg(Y\times Z))$.
Write $R_i=\sum_{\a}R_{i,\a}\ot \ell^{\a}$ with $R_{i,\a}\in I_2(Y)$, so
$$
P(v\ot z,v\ot z, w\ot y)   
=\frac 32\sum R_{i,\a}(v,v)  \ell^\a(w)T_i(z,z,y) 
$$
which  is zero because $R_{i,\a}\in I_2(Y)$ for all $i,\a$.
  Similarly $\pi_S(K_S(A)\ot S_1(Z))\subset I_3(\s(Seg(Y\times Z))$.

Now say $P=\pi_S(\sum R_i\ot T_i)\in (K_S(A)\ot K_S(B))\cap I_3(\s(Seg(Y\times Z))$.
Without loss of generality we assume   the $R_i$ are linearly independent modulo $S_1(Y)$
and   the $T_i$ are linearly independent modulo $S_1(Z)$. 
Fix $(y,z)\in U$, we obtain a linear equation
$$
\sum_i R_i(v,v,w)c_{i,y,z} =0
$$
where $c_{i,y,z}=T_i(y,y,z)$ and if $y,z$ are chosen generically all the coefficients are nonzero because we
are working mod $S_1(Z)$. Note that the index range for
$i$ is at most from $1$ to $\tmin\{\tdim K_S(A),\tdim K_S(B)\}$. 
We will show each $R_i(v,v,w)$   must be zero for all $v,w\in\hat Y$.
Since $\hat Y$ spans $A$ and the expression is linear in $w$, we can have this hold for all $v\in \hat Y$ and
$w\in A$, but this in turn implies that each $R_i\in S_1(Y)$.  

To obtain the desired vanishing, fix $v,w$ and consider the $R_i(v,v,w)=r_i$ as constants. We have
an equation
$$
\sum_i r_iT_i(y,y,z)=0 \ \ \forall y,z\in \hat Z
$$
As remarked above, since $Z$ is linearly non-degenerate, we may choose $\tdim B$ elements
$z_s\in \hat Z$ that give a basis of $B^*$. Similarly, we may choose  $\binom{\tdim B+1}2-\tdim I_2(Z)$
elements $y_t\in \hat Z$ such that the vectors  $y_t^2$ span $ I_2(Z)\upperp\subset S^2B^*$.
Thus the vectors $y_t^2\ot z_s$ give a basis of $ I_2(Z)\upperp\ot B^*$. Thus the pairing with elements  of
$K_S(B)/S_1(Z)$ is perfect, which implies that the matrix given by pairing the $y_t^2\ot z_s$ with the $T_i$
has a one-sided inverse, so we have enough independent equations to force all the $r_i$ to vanish.

\smallskip

The argument for the $S^3A\ot S^3B$ factor is similar, but easier, as there is no
need to symmetrize. Write 
$P= \sum R_i\ot T_i $ with $  R_i\in S^3A$ and $  T_i\in    S^3B$.
$$
P(v\ot z,v\ot z, w\ot y)   
=\sum R_{i}(v,v,w)   T_i(z,z,y).
$$
Now
$R_i\in I_3(\s(Y))$ iff $R_i(v,v,w)=0$ for all $v,w\in \hat Y$ and one concludes as above.
\end{proof}

\subsection*{proof of Theorem \ref{dcomthm}}
We first observe that without loss of
generality, we may assume $\tdim A=2$.
This is because, continuing the notation of \S\ref{desingsect}, 
$\s(Seg(\BP A^*\times Y))\subset \s(Seg(\BP A^*\times \BP W^*))$
and the ideal of the latter already contains
all partitions of length greater than two and
is generated by
$\La 3A \ot \La 3 W $.
Thus the only new modules of generators that occur  when
$\tdim A>2$ are the components of  $\La 3A \ot \La 3 W $.

The modules $I_2(\s(Y))$   are all trivial modules except for
$V_{\o_6}^{A_n}$ for $n\geq 6$ (and $\La 2 B\ot \La 2 C$ for
$Y=Seg(\BP B^*\times \BP C^*)$.

We now show there are no
new generators in degrees greater than three.

\subsection{Case $Y=G(2,B)$} We give two proofs of this case, the first gives more
information, the second is uniform with the case $Y=\BO\pp 2$.
  
\subsubsection*{First proof} 
We need to study the exactness of the middle step of
$$ 
(S_{11}A\ot S_{11}(\La 2B))\ot   
S^{r-2}(A\ot \La 2 B )\mid_{|\pi |\leq 4}   \ra   
(A\ot \La 2 B)\ot  S^{r-1}(A\ot\La 2 B)\mid_{|\pi |\leq 4}   \ra   
 S^{r }(A\ot\La 2 B)\mid_{|\pi |\leq 4}
$$
Here by $\mid_{|\pi |\leq 4}$, we mean the components
of $S_{a,b}(\La 2 B)$ that occur
in the decomposition of  $S^p(A\ot \La 2B)=\oplus_{a+b=p}S_{a,b}A
\ot S_{a,b}(S_{1,1}B)$ that as partitions $S_{\pi}B$ have
length at most four.  

Since the   partitions $S_{\pi}B$ occurring
in the middle  entry can have length at most six,
it is sufficient to solve the problem for the
case $n\leq 6$.  The decomposition of $S_{a,b}(\La 2B)$ is not
known in closed form,   however at this point
we could rely on a computer to compute $H^{d-1}(\BP A^*\times G(4,6), \La d\xi)$.
We instead use an induction argument that is computer free.

Let $B$ have dimension $n$ and consider the rank variety
$$
R_{n-1}(A^*\ot \La 2 B^*):=
\{ T\in A^*\ot \La 2 B^* \mid \exists U,\
\tdim U=n-1,\ T\in A^*\ot \La 2 U\}.
$$
 The
secant variety $\hat\s(\BP A^*\times G(2,B^*))$ coincides
with  
$R_4(A^*\ot \La 2 B^*)$.
We determine the ideal of $R_n$ via that of $R_{n-1}$
which will render the bundles $\xi$ that we use irreducible.

Over $\cB:=A^*\ot G(n-1,B^*)$ we have
the bundle with fiber $A^*\ot \La 2\cS$ which provides
a desingularization of $R_{n-1}(A^*\ot \La 2 B^*)$.
Our corresponding bundles are
$\eta=A\ot \La 2\cS^*$ and $\xi=A\ot (\La 2 B/\La 2 \cS)
=A\ot \cS^*\ot \cQ^*$.
Note that $\trank \xi=2(n-1)$.

We have 
$$
\La d\xi=\bigoplus_{\{\pi= (a,b)\mid 2a+b=d,\ a+b\leq n-1\}}
S_{\pi}A\ot S_{\pi'}\cQ^*\ot S_d\cS^*
$$
where $\pi'$ denotes the conjugate partition 
to $\pi$, so if $\pi=(a,b)$, then
$\pi'=(2^a,1^b)$.

We apply the Bott algorithm
(see e.g.  \cite{weyman}, \S 4.1.5)   to
the weight
$\mu=\o_a+\o_{a+b}-(2a+b)\o_{n}$
The only
potential way
to have non-zero cohomology is at step $(n-1)-(a+b)$, or
at step $(n-1)-a$ or step $n-1$.

To  obtain nonzero  $H^{(n-1)-(a+b)}(\cB,\La{n-(a+b)}\xi)$,  
$-(2a+b)+(n-1-(a+b))$ must be negative and
$-(2a+b)+(n-1-(a+b))+2$ must be non-negative,
so we must have $3a+2b=n-1$.
Write   $i=2a+b$, so
$H^{i}(\cB, \La{i+1}\xi)=S_{n-i,2i-n}A\ot S_{2^{2i-n},1^{2n-2i}}B$.

To obtain nonzero  $H^{(n-1)-a}(\cB,\La{n-a}\xi)$,  
$-(2a+b)+(n-1-a)$ must be negative and
$-(2a+b)+(n-1-a)+2$ must be non-negative,
implying $3a+b=n+1$, contradicting $a+b\leq n-1$.

To obtain nonzero  $H^{ n-1 }(\cB,\La n\xi)$, we
would have to have
$-(2a+b)+(n-1)<-1$   and
$-(2a+b)+(n-1)+2$ non-negative,
implying $2a+b=n+2$ contradicting $a+b\leq n-1$.

In summary:
\begin{proposition}  
Let $B,A$ be vector spaces respectively  of dimensions $n,2$   and 
consider the rank variety
$$
R_{n-1}(A^*\ot \La 2 B^*):=
\{ T\in A^*\ot \La 2 B^* \mid \exists U,\
\tdim U=n-1,\ T\in A^*\ot \La 2U\}.
$$
Then the ideal of $R_{n-1}$ is generated in degrees
$\lceil \frac n2\rceil \leq d\leq \lfloor \frac{2n}3\rfloor$  
by the modules
$$
S_{n-d,2d-n}A\ot S_{2^{2d-n},1^{2n-2d}}B
\subset S^d(A\ot \La 2B)
$$
\end{proposition}

To prove Theorem \ref{dcomthm}, by the discussion above we need to examine the
cases $n=5,6$. When $n=5$ we only have $d=3$
and the module
$S_{2,1}A\ot S_{2,1^4}B$ generates the ideal
of $R_4(A^*\ot \La 2 B^*)$ and also of $\s(\BP A\times G(2,5))$.
When $n=6$ we have $d=3,4$ and the modules
$S_{3}A\ot S_{1^6}B$
in degree $3$ and
$S_{2,2}A\ot S_{2^2,1^4}B$ in
degree $4$ generate the ideal of $R_5(A^*\ot \La 2 B^*)$,
and the ideal of 
$\s(Seg(\BP A\times G(2,6)))$ is generated by these
and the representation
$S_{2,1}A\ot S_{2,1^4}B$ that already occurs for $dim\ A=5$. However,
$S_{2,2}A\ot S_{2^2,1^4}B$ is already in the ideal
generated by 
$S_{2,1}A\ot S_{2,1^4}B$.

\subsection{Proof   of cases $Y=G(2,K^6)$ and $Y=\BO\pp 2$}
Write $\cB=G/P_{i_0}$. 
Following the conventions of \cite{bour}, $i_0=4,6$ respectively.
We write the Levi factor of $\fp$ as $\fg_0=\ff + \langle Z_{i_0}\rangle$, where
$\ff$ is semi-simple (respectively $\fa_3+\fa_1$ and $\fd_5$) and
$\langle Z_{i_0}\rangle $ is the center of $\fg_0$.

We will obtain the result by computing $H^i(\cB,\La{i+1}\xi)$ via
$H^i(\cB,\La{i+1}gr(\xi))$ and applying a result of Ottaviani and Rubei.
 
Here  
\begin{align*}
 gr(\xi)&=A\ot (\cS^*\ot \cQ^*\op \La 2 Q^*)=A\ot (E_{\o_1-\o_4 +\o_5}\op E_{-\o_4}) \ 
{\rm for\ }Y=G(2,6)\\
&= A\ot (E_{\o_2-\o_6}\op E_{-\o_6}) \ {\rm for\ } Y=\BO\pp 2
\end{align*}
where $E_{\l}$ denotes the irreducible bundle corresponding to the $\fg_0$-module
of highest weight $\l$.

We first compute the decomposition of the
exterior powers of the  $\fg_0$-module giving rise to $gr(\xi)$ as an
$\ff$  module and then compute the action of $Z_{i_0}$ to determine
  the coefficient  on  $\o_{i_0}$   for each irreducible $\ff$-module appearing.

The $\ff$-module decomposition is straightforward with the aid of LiE \cite{LiE},
keeping in mind that  $\tdim A=2$:
\begin{align}
 \La k (A\ot (U\op K))
=&\op_{a+b=k} S_{a,b}A\ot S_{2^a,1^b} (U\op K)\\
=&\op_{a+b=k} S_{a,b}A\ot (S_{2^a,1^b}U\op  
  S_{2^a,1^{b-1}}U\op S_{2^{a-1},1^{b+1}}U
\oplus S_{2^{a-1}, 1^b}U) 
\end{align}

One then uses LiE to decompose these GL(U)-modules as $\ff$-modules.
Next to determine the weight on  the marked node (i.e., the coefficient
of $\o_{i_0}$), one uses
the grading element $Z_{i_0}\in \ft$ which has the property
that $Z_{i_0}(\a_j)=\d_{i_0,j}$. 
Thus if $\l=\sum_{i\neq i_0}\l^i\o_i$ is an irreducible $\ff$-module
appearing in $W_{\mu}^{\ot k}$, where the $\o_i$ are fundamental weights of $\fg$,
to find the coefficient of $\o_{i_0}$ of the $\fg_0$-module, one calculates
$$
\sum\l^j(c\inv)_{i_0,j}  = k Z_{i_0}(\mu)
$$
where $(c\inv)$ denotes the inverse of the Cartan matrix.
In both our cases $Z_{i_0}(\mu)=-\frac 13$.

Now one calculates  $H^j(\cB,\La p gr(\xi))$.
In practice we first calculated $H^{p-1}(\cB, \La{p} gr(\xi))$, and only
if this was nonzero did we calculate the other $H^j(\cB, \La p gr(\xi))$.

We got no relevant cohomology except in degrees one and two. In
both cases the only modules that appeared in degree two were the
cubic generators of the ideal plus 
$S_{21}A\otimes \BC$, which is cancelled by
its appearance in $H^1$, and it is the unique module
appearing in $H^1$.  By Proposition 6.7 in \cite{OR}  
cancellation occurs in the spectral sequence   for the cohomology of $\La 3\xi$.
The program we used for this calculation is publicly available
at www.math.tamu.edu/$\sim$robles.

\subsection{$Y=\BS_5$}
We need to calculate $H^i(\La{i-1}\xi)$ with $\xi=\BC^2\ot E$ where
$E$ is the vector bundle determined by the $P$-module with highest
weight $\l=[-1,0,0,1,0]$. This calculation is similar to the above,
but significantly easier because $\xi$ is irreducible.

\section{The coordinate ring of $\s (X)$}\label{coordringsect}

The following proposition is   due to F. Zak (\cite{zak}, p. 51):

\begin{proposition} Let $X=G/P\subset \BP V$ be a homogeneously
 embedded homogeneous variety. Let $\l$ denote the highest weight and $\mu$ denote the lowest 
weight of $V$, and let $v_{\l},v_{\mu}$ be corresponding weight vectors.
Then   $\s(X)=\overline{G.[v_{\l}+v_{\mu}]}$, where the closure is the Zariski closure.
\end{proposition}
\begin{proof}
$$\fg.(v_{\l}+v_{\mu})=
(\fg_+ + \fg_0 + \fg_{-}).(v_{\l}+v_{\mu})
=\fg_+.v_{\l} + \fg_{-}.v_{\mu} =\hat T_{[v_{\l}+v_{\mu}]}\s(G/P)
$$
where the last equality is Terracini's lemma. This proves the result
in the case $\fg_+.v_{\l}\cap \fg_{-}.v_{\mu}=0$, i.e. the secant variety
is non-degenerate, and the result is easy to verify in the degenerate cases 
as well. (The degenerate cases are all rank at most 2 CHSS and the adjoint varieties $G/P_{\tilde\a}\subset\BP \fg$
where $\tilde\a$ is the largest root.)
\end{proof}
 
We work with the affine variety $ \hat \s (X)\subset V$.
\begin{proposition} Notations as above.
 Let $H=Stab( v_{\l}+v_{\mu} )$.
Then 
$$K[ \hat \s (X)] \subseteq K [ G.(v_{\l}+v_{\mu})] = K[G/H]=
K [G]^H = \oplus_{\l\in\Lambda^+_G}V_{\l}\ot V^{*H}_{\l}. 
$$
In particular,  the irreducible $G$-module
$V_{\l}$ occurs  in $K [G]^H$ with multiplicity equal
to the number of $H$-fixed points in $V_{\l}^*$.
\end{proposition}
\begin{proof}
Since $G/H\subseteq \hat \s(X)$ we obtain an inclusion
$K[ \hat \s (X)]\subseteq K[G/H]$ by restricting functions on $\hat \s(X)$.
By \cite{[K]}, Theorem 3, Chapter II, section 3, the coordinate ring
of $G$ has a left-right decomposition (as a $(G-G)$-bimodule)
$$K[G]= \oplus_{\{\l\in \Lambda^+_G\} } V_{\l}\ot V_{\l}^*, 
$$
    taking $H$-invariants we get
$$K[G]^H= \oplus_{\{\l\in \Lambda^+_G\} } V_{\l}\ot (V_{\l}^*)^H.$$ 
\end{proof}

Note that if   $G/H\subset \BP V$, then   $K[G]^H$ is equipped with a grading
(that depends on the embedding). We do not know
of any way to recover the grading from this description in general, in fact the
same $V_{\l}$ may appear in several different degrees. Fortunately, when $G=GL_n$,
the degree is recoverable. 
 
\begin{example} $X=G(k,W)$, with $k>2$.
Here without loss of
generality we may take
$\tdim W\ge 2k$. Indeed, if  $dim W< 2k$, we may pass to the dual Grassmannian $G(\tdim W-k, W^*)$. If $dim(W)>2k$, $\hat \s(X)$ is contained in the subspace variety $R_{2k}(\La k W)$ of tensors that can be written using $\le 2k$ basis vectors.  
Let $W'\subset W$ be a $2k$ dimensional subspace.
Consider the subgroup $H'\subset H$,

$$H'=\lbrace \phi\in SL_K (W)\ |\ \phi |_{W'}=Id_{W'}\}.$$

The quotient $SL_K (W)/H'$ can be identified with the variety $Hom_K^{inj}(W' ,W)$ of injective linear maps from $W'$ to $W$. Since the complement of  $Hom_K^{inj}(W' ,W)$ in $Hom_K (W' ,W)$ has codimension $\ge 2$ every regular function on  $Hom_K^{inj}(W' ,W)$ extends to $Hom_K (W' ,W)$. This means that we have the equalities

$$K[R_{2k}(\La k W)]=K[SL_K (W)/H' ]=\oplus_{\l\in\Lambda^+_{SL_K(W)}} (S_\l W)^*\otimes S_\l W^{H'}=\oplus_{\l} (S_\l W)^*\ot S_\l W'.$$
Here  in the last equality we may view $\l$ as a partition. Note that the last equality states that
the module $(S_{\l}W)^*$ appears with multiplicity $\tdim S_{\l}W'$.

This reduces the calculation of $(S_\l W)^H$ to the case $W=W'$.
Assuming now that $dim (W)=2k$, with basis $e_1\hd e_{2k}$,
we may take $ v_{\l}+v_{\mu}= e_1\ww \cdots \ww e_k
+e_{k+1}\ww\cdots\ww e_{2k}$.
Then
$$
H= \left\{ \begin{pmatrix} A&0\\ 0&B\end{pmatrix}
\mid det(A)=det(B)=1\right\}  
$$

Considering $G(k,W)=SL(W)/P_k$, 
it is  clear that no fundamental representation other than
$W_{\o_k}=\La kW$ has an $H$-fixed vector and in
$W_{\o_k}$ there is a $2$-dimensional subspace of such spanned by
$e_1\ww\cdots \ww e_k$ and $e_{k+1}\ww\cdots\ww e_{2k}$. The corresponding two copies of $\La k W$ generate the ring of invariants in the following sense.
We claim:
$$
  K [SL(W)]^H=  \oplus_{r,s\ge 0} S_{(r^k)}W\otimes S_{(s^k)}W.
$$
Two generating fundamental representations are in bidegrees $(1,0)$ and $(0,1)$.
If we   work instead with $GL(W)$, since $H$ acts trivially on the determinant, we get
$$
  K [GL(W)]^H=  \oplus_{r,s\ge 0, m\in \BZ} S_{(r^k)}W\otimes S_{(s^k)}W\ot (\det)^m.
$$
This second description has the advantage that when we consider the embedded
space $G/H\subset \BP \La k W$ we can determine the degree these modules appear in.

To see this, write $E=\langle e_1\hd e_k\rangle, F=\langle e_{k+1}\hd e_{2k}\rangle$,
we want to see how many instances of the trivial representation of
$SL(E)\times SL(F)$ occurs in the irreducible $SL(W)$ module $S_{\l}W$.
Now, since $W=E\op F$ 
$$S_\lambda (W)=\oplus_\mu S_\mu E\otimes S_{\lambda/\mu}F
=\oplus_{\mu ,\nu}c^\l_{\mu ,\nu} S_\mu E\ot S_\nu F
$$

we have
$${\rm dim}\ (S_\l W^* )^H = \sum_{r,s\ge 0} c^{\l}_{(r^k), (s^k)} .$$
This gives

$$K[SL(W)]^H =\oplus_{r,s\ge 0} S_{(r^k)}W\ot S_{(s^k)}W .$$
This means this ring is the homomorphic image of the symmetric algebra on two copies of $\La kW$ corresponding to components in bidegrees $(1,0)$ and $(0,1)$ which is our claim.

\end{example}

\begin{example} $X=Seg(\BP A_1\ot \cdots\ot \BP A_k)=\Pi GL(A_i)/P$, $\tdim A_i=2$.
Here we may take $ v_{\l}+v_{\mu}= e_1\ot \cdots \ot e_k
+f_1\ot\cdots\ot f_k$  with $e_j,f_j$ a basis of $A_j$.
Then
$$
H= \left\{ \Pi_j\begin{pmatrix} s_j&0\\  0&t_j  \end{pmatrix}
\mid s_1\cdots s_k=t_1\cdots t_k=1 \right\}
$$ 
Now consider the action of $H$ on
$$S_{a_1,b_1}A_1\ot \cdots\ot S_{a_k,b_k}A_k=
(det A_1)^{b_1}S^{a_1-b_1}A_1\ot \cdots \ot (det A_k)^{b_k}S^{a_k-b_k}A_k
$$
The weight vectors are
$$e_1^{i_1+b_1}f_1^{a_1 -i_1}\ot \cdots \ot 
e_k^{i_k}f_k^{a_k- i_k},  \ 0\leq i_{j}\leq a_{j} -b_{j}
$$
which is acted on by $H$ by
$$(s_1^{i_1+b_1}  \cdots   
s_k^{i_k+b_k})(t_1^{a_1-b_1-i_1}\cdots t_k^{a_k-b_k-i_k} ) 
$$
so we need
$ 
i_1+b_1=i_2+b_2=\cdots = i_k+b_k$
and
$a_1-b_1-i_1=a_2-b_2-i_2=\cdots = a_k-b_k-i_k$.
This means the vectors $e_j^{i_j+b_j}f_j^{a_j-i_j}$ have to be all of the same weight, for $j=1,\ldots ,k$.

This means for a module
$$S_{a_1 ,b_1}A_1\ot \cdots\ot S_{a_k ,b_k}A_k,$$
the dimension of the subspace of $H$-invariant vectors is ${\rm min}\ a_j -{\rm max}\ b_j +1$.

\end{example}

\end{document}